\documentclass[11pt,a4]{article}
 \usepackage{amsmath,amssymb,latexsym,hyperref,xcolor,amsthm,authblk}
 \def\le{\leqslant}
 \def\ge{\geqslant}
 \def\Z{\mathbb Z}
 \def\N{\mathbb N}

 \def\epsilon{\varepsilon}
 \def\phi{\varphi}
  \newtheorem{defi}{Definition}
  \newcommand{\bd}{\begin{defi}}      
  \newcommand{\ed}{\end{defi}}

  \newtheorem{lemm}[defi]{Lemma}  
  \newcommand{\bl}{\begin{lemm}}
  \newcommand{\el}{\end{lemm}} 

  \newtheorem{theo}[defi]{Theorem}
  \newcommand{\bt}{\begin{theo}}
  \newcommand{\et}{\end{theo}}

  \newtheorem{cor}[defi]{Corollary}
  \newcommand{\bc}{\begin{cor}}
  \newcommand{\ec}{\end{cor}}

  \newtheorem{pro}[defi]{Proposition}
  \newcommand{\bp}{\begin{pro}}
  \newcommand{\ep}{\end{pro}}
  \theoremstyle{definition}
    \newtheorem{rmk}[defi]{Remark}
  \newcommand{\brmk}{\begin{rmk}}
  \newcommand{\ermk}{\end{rmk}}
  
  \newtheorem*{theo*}{Theorem}
  
      \newtheorem{ex}[defi]{Example}
  \newcommand{\bex}{\begin{ex}}
  \newcommand{\eex}{\end{ex}}

  \makeatletter
  \def\proof{\@ifnextchar[\opargproof{\opargproof[\bf Proof \hfil\\ ]}}
  \def\opargproof[#1]{\par\noindent {\bf #1 }}
  \def\endproof{{\unskip\nobreak\hfil\penalty50\hskip8mm\hbox{}
  \nobreak\hfil
  \(\Box\)\parfillskip=0mm \par\vspace{3mm}}}
  \makeatother
  
\begin{document}
\title{\LARGE\bf Polyfunctions over Commutative Rings}
\author[1]{Ernst Specker}
\author[2]{Norbert Hungerb\"uhler}
\author[3]{Micha Wasem}
\affil[1]{\em Dedicated to the memory of the first author\bigskip}
\affil[2]{Department of Mathematics, ETH Z\"urich, R\"amistrasse 101, 8092 Z\"urich, Switzerland\bigskip}
\affil[3]{HTA Freiburg, HES-SO University of Applied Sciences and Arts Western Switzerland, P\'erolles 80, 1700 Freiburg, Switzerland}

\maketitle
 \begin{abstract}
\textcolor{black}{A function $f:R\to R$, where $R$ is a commutative ring with unit element, is called \emph{polyfunction} if it admits a polynomial representative $p\in R[x]$.
Based on this notion we introduce ring invariants which associate to $R$ the numbers $s(R)$ and $s(R';R)$, 
where $R'$ is the subring generated by $1$.  For the ring $R=\mathbb Z/n\mathbb Z$
the invariant $s(R)$ coincides with the number theoretic \emph{Smarandache} or {\em Kempner function} $s(n)$. 
If every function in a ring $R$ is a polyfunction, then $R$ is a finite field according to the 
R\'edei-Szele theorem, and it holds that $s(R)=|R|$. However, the condition $s(R)=|R|$ 
does not imply that every function $f:R\to R$ is a polyfunction. We classify all finite 
commutative rings $R$ with unit element which satisfy $s(R)=|R|$.  For  infinite rings $R$, 
we obtain a bound on the cardinality of the subring $R'$ and for $s(R';R)$ in terms of $s(R)$.
In particular we show that $|R'|\leqslant s(R)!$. We also give two new proofs for the 
R\'edei-Szele theorem which are based on our results.}
 \end{abstract}
 
\section{Introduction}

\textcolor{black}{For a commutative ring $R$ with unit element, a function $f:R\to R$ is said to be a \emph{polyfunction} if there exists a polynomial $p\in R[x]$ such that $f(x)=p(x)$ for all $x\in R$ (see~\cite{redei,mullen}, and also \cite{bhargava,chen} for a discussion on polyfunctions from $\Z_m\to\Z_n$). The set of polyfunctions over $R$ equipped with pointwise addition and multiplication forms a subring
$$G(R):=\{f:R\to R, \exists p\in R[x]\,\, \forall x\in R\implies p(x)=f(x)\}$$
of $R^R$ and will be called the \emph{ring of polyfunctions} over $R$. The polynomials in $R[x]$ which represent the zero element in $G(R)$ are called \emph{null-polynomials} (see~\cite{singmaster}).
If $S$ is a subring of $R$, then
$$G(S;R):=\{f:R\to R, \exists p\in S[x] \,\,\forall x\in R\implies p(x)=f(x)\},$$
is a natural subring of $G(R)$. In particular, the subring $R'$ generated by the unit element 1
in $R$ gives rise to the {\em integer polyfunctions\/}
$G(R';R)$.
Instead of restricting the ring of allowed coefficients as in
the construction for $G(S;R)$,
one obtains other rings of polyfunctions by restricting the domain: The ring $$
\{f:S\to R, \exists p\in R[x]\,\,\forall x\in S\implies p(x)=f(x)\}
$$ e.g.\ contains $G(R)$ as a subring.}


\textcolor{black}{If $S$ is a subring of $R$}, a characteristic number connected to $S$ and $R$ is the minimal degree
$m$ such that the function $x\mapsto x^{m}$ can be represented by a 
polynomial in $S[x]$ of degree \textcolor{black}{strictly} smaller than $m$. Then, in particular,
every function in $G(S;R)$ has a polynomial representative
of degree strictly less than $m$.
We set
$$s(S;R):=\min\{m\in \N, \exists p\in S[x], 
\operatorname{deg}(p)<m, \forall x\in R \implies p(x)=x^{m}\}$$
and $s(R):=s(R;R)$ for brevity. We set $s(S;R):=\infty$ if no  function $x\mapsto x^m$ can be 
represented by a polynomial of degree strictly smaller than $m$.

Trivially, we have $s(S;R)\ge s(T;R)\ge s(R)$ whenever $S\subset T$ are
subrings of $R$. On the other hand, 
we will see in Section~\ref{generali}, 
that $s(R';R)<\infty$ is bounded in terms of $s(R)$
if $s(R)<\infty$.

Clearly, if two rings $R_1,R_2$ are isomorphic, then $s(R_1) = s(R_2)$ and $s(R_1',R_1)= s(R_2',R_2)$.
In other words, $R\mapsto s(R)$ and $R\mapsto s(R',R)$ are ring invariants.

\textcolor{black}{The function $s$,
which associates to a given ring $R$ the number $s(R)\color{black}{\in \N\cup\{\infty\}}$ has been introduced in~\cite{artikel1} and is 
called {\em Smarandache function}. This naming stems from the fact, that for all $2\leqslant n\in\mathbb N$, the map} $n\mapsto s(\Z/n\Z)$
coincides with the well-known number theoretic 
Smarandache or Kempner function $s$ (see~\cite[Theorem 2]{artikel1}) defined by
\begin{equation}\label{standardsmarandache}
s(n) := \min\{k\in\N,n\mid k!\}
\end{equation}
(see Lucas~\cite{lucas}, Neuberg~\cite{neuberg} and Kempner~\cite{kempner}).
In fact, Legendre has already studied aspects of the function $s(n)$: In~\cite{legendre}
he showed that if $n=p^\mu$ for some prime $p$ and $1\leqslant \mu\in\N$, then $s(n)$ verifies
$$
s(n)=\mu(p-1)+a_0+a_1+\ldots+a_k,
$$
where the numbers $a_i$ are the digits of $s(n)$ in base $p$. i.e.\ $s(n) = a_kp^k+\ldots + a_0$ and $0\leqslant a_i < p$.
We refer to Dickson~\cite[p.\ 263--265]{dickson} for the history of the function $s(n)$.

\textcolor{black}
{In a finite field $F$, every function is a polyfunction as a polynomial respresentative of a function $f:F\to F$ is, e.g., given by the Lagrange interpolation polynomial for $f$. This representation property characterizes finite fields among commutative rings with unit element (see~\cite{redei2}):}

\bt[R\'edei, Szele]\label{iff}
If $R$ is a commutative ring with unit element then $R$ is a finite field if and only if every function $f:R\to R$ can be represented by a polynomial in $R[x]$.
\et
\textcolor{black}{We will include two short alternative proofs of this theorem in Section~\ref{alternative}. For finite fields $F$, one has $s(F)=|F|$, so in view of Theorem~\ref{iff}, it is natural to ask what can be said about commutative rings $R$ with unit element for which $s(R)=|R|$ holds true. \textcolor{black}{Note that if $R$ is a finite ring, it trivially holds that $s(R)\leqslant |R|$, as the polynomial
$$
p(x) = \prod_{y\in R}(x-y)
$$
is a normed null-polynomial of degree $|R|$.}
}

\textcolor{black}{The following theorem (which will be restated below for the reader's convenience as Theorem~\ref{oppo}), 
answers the above question and classifies all finite commutative rings $R$ with unit element that satisfy $s(R)=|R|$:}
\begin{theo*}\textit{Let $R$ be a finite commutative ring with unit element.
Then, $s(R)=|R|$ holds if and only if $R$ is one of the following:
\begin{itemize}
\item[(a)] $R$ is a finite field, or
\item[(b)] $R$ is $\Z_4$, or
\item[(c)] $R$ is the ring $\rho$ with four elements $\{0,1,a,1+a\}$
with $1+1=0$ and $a^2=0$.
\end{itemize}}
\end{theo*}

\textcolor{black}{{\bf Remarks:}
\begin{enumerate}
\item The ring $\rho$ is not a field since it has zero divisors, and since it is of characteristic $2$, it is not isomorphic to $\Z_4$.
\item Observe the similarity between this result and the fact that for $n\geqslant 2$, the usual Smarandache function satisfies $s(n)=n$ if and only if $n$ is prime or $n=4$.
\end{enumerate}}
\textcolor{black}{Section~\ref{generalf} is devoted to the proof of this theorem.} In Section~\ref{generali} we discuss infinite rings and show that for an infinite commutative ring $R$ with unit element and $s(R)<\infty$, we obtain an upper bound for $|R'|$ and for $s(R'; R)$ \textcolor{black}{in terms of $s(R)$}, where $R'$ is the subring of $R$ generated by $1$. 
\textcolor{black}{Finally, in Section~\ref{alternative}, we give two proofs of Theorem~\ref{iff} -- a direct one and one that is based on Theorem~\ref{oppo}.}

\textcolor{black}{Throughout the article, $n\ge 2$ will denote a natural number, and $\Z_{n}=
\Z/n\Z$\, is the ring of integers modulo $n$, and we write $a\mid b$ if $b$ is an integer multiple of $a$. }

\section{Polyfunctions over Finite Rings}\label{generalf}

Theorem~\ref{iff} answered the question, when
a ring $R$ has the property, that every function $f:R\to R$ 
can be represented by a polynomial in $R[x]$. For finite
rings a necessary (but not sufficient) condition 
for this property to hold is 
\begin{equation}\label{imp}
s(R)=|R|,
\end{equation}
(see Theorem~\ref{oppo} below). In this section, we want to address the question for which
finite rings, equation~(\ref{imp}) holds. The first step 
to answer this, is the following proposition:

\bp\label{vergessen}
If $R$ is a commutative ring with unit element
and with zero divisors then either
\begin{itemize}
\item[(a)] there exist $a,b\in R\setminus\{0\}$ with $a\neq b$
and $ab=0$, or
\item[(b)] $R$ is $\Z_4$, or
\item[(c)] $R$ is the ring $\rho$ with four elements $\{0,1,a,1+a\}$
with $1+1=0$ and $a^2=0$.
\end{itemize}
\ep
\proof
Let us assume that in $R$ the implication holds:
if $u,v\in R\setminus\{0\}$ and $uv=0$ then it follows $u=v$.
Let $a\in R\setminus\{0\}$ be a zero divisor: $a^2=0$. 
Thus, if $x$ is an element in $R$ with $ax=0$, we have either
$x=0$ or $x=a$. Notice that for all $u\in R$ we have
$$ a(au)=0$$
and hence  for all $u\in R$
$$au=0\text{ or $a(u-1)=0$.}$$
Hence, we have only the four cases $u=0$ or $u=a$ or $u=1$ or $u=1+a$.
If $1+1=a$, then $R=\Z_4$, if $1+1=0$, then $R$ is the ring $\rho$ in (c).
\endproof

We can now prove the main result of this section:
\bt\label{oppo}
Let $R$ be a finite commutative ring with unit element.
Then, $s(R)=|R|$ holds if and only if $R$ is one of the following:
\begin{itemize}
\item[(a)] $R$ is a finite field, or
\item[(b)] $R$ is $\Z_4$, or
\item[(c)] $R$ is the ring $\rho$ with four elements $\{0,1,a,1+a\}$
with $1+1=0$ and $a^2=0$.
\end{itemize}
\et

\proof
If $R$ is not a field and not $\Z_4$ and not the ring $\rho$, then, according
to Proposition~\ref{vergessen}, $R$ is a ring with
$a,b\in R\setminus\{0\}$ such that $ab=0$ and with $a\neq b$. Then 
$$
(x-a)(x-b)\prod_{z\in R\setminus\{a,b,0\}}(x-z)
$$
is a normed null-polynomial of degree $|R|-1$. Therefore
$s(R)<|R|$.

To prove the opposite direction, we go through the three cases:

(a) If $R$ is a field, then a polynomial of degree $n$ has at most
$n$ roots. Hence, $s(R)=|R|$.

(b) If $R$ is $\Z_4$, then \textcolor{black}{(by \cite[Theorem 2]{artikel1})
$s(\Z_4)=s(4)=4=|\Z_4|$.}


(c) If $R$ is the ring $\rho$ with elements $\{0,1,a,1+a\}$
and with $1+1=0$ and $a^2=0$, we have to prove that
$s(R)=4$. Assume by contradiction, that $p(x)\in R[x]$ is 
a normed null-polynomial of degree  3.
Since $p(0)=p(1)=0$, $p(x)$ must be of the form
$$ p(x)=x(x+1)(\xi+x).$$
From $p(a)=0$, it follows that 
$a\xi=0$ and from $p(a+1)=0$ it subsequently follows that
$a=0$ which is a contradiction.
\endproof

\section{Infinite Rings}\label{generali}
In this section $R$ is a commutative ring with unit element and
$R'$ the subring of $R$ which is generated by $1$. \textcolor{black}{We will need the following lemma, which is a corollary of \cite[Lemma 4, p.4]{artikel1}:
\bl\label{kombi}
For all $k,n\in \N\cup\{0\}, k\leqslant n$ we have
$$\sum_{j=0}^n (-1)^{n-j}\binom nj j^k = \delta_{kn}n!$$
(with the convention $0^0:=1$).
\el}
\bp\label{fini}
If $s(R)<\infty$ then $R'$ is a finite ring and
$|R'|\big| s(R) ! $.
\ep
{\bf Remark:} We notice, that $s(R)<\infty$ may hold even
if $R$ is an infinite ring. As an example consider the ring
$$
R=\Z_2[x_1,x_2,\ldots]/\{x_1^2,x_2^2,\ldots\}
$$
in which all $u\in R$ satisfy the relation $u^4=u^2$.
On the other hand, if $R$ is finite, we trivially
have $s(R)\le |R|$.

{\bf Proof of Proposition~\ref{fini}}\\
By assumption, for $n=s(R)$ there exist coefficients $a_i\in R$, $i\in
\{0,1,\ldots n-1\}$, such that for all $u\in R$ we have 
\begin{equation}\label{assu}
u^n-\sum_{i=0}^{n-1}a_i u^i=0.
\end{equation}
We denote 
$$\underbrace{1+1+\ldots+1}_{{\text{$m$ times}}}\in R'$$ 
by $\bar m$.
Then, by Lemma~\ref{kombi}, we have for $k\le n$
\begin{equation}\label{qw}
\sum_{j=0}^n \overline{(-1)^{n-j}\binom nj j^k} = \overline{\delta_{kn}n ! }
\end{equation}
Hence, it follows from~(\ref{assu}) that
\begin{eqnarray*}
0 &=& \sum_{j=0}^n \overline{(-1)^{n-j}\binom nj}\left(\bar j^n-
\sum_{i=0}^{n-1}a_i\bar j^i\right)=\\
&=&\sum_{j=0}^n \overline{(-1)^{n-j}\binom njj^n}-
\sum_{i=0}^{n-1}a_i\sum_{j=0}^n \overline{(-1)^{n-j}\binom njj^i}=\overline{n ! }
\end{eqnarray*}
where the last equality follows from~(\ref{qw}).
\endproof
{\bf Remark:} As the example $R=\Z_{n!}$ shows, the estimate
on the size of $R'$ emerging from
Proposition~\ref{fini}, $|R'|\le s(R)! $,  cannot be improved in general.
\bl\label{mul}
If $n:=s(R)<\infty$ then there exists a bound $\color{black}{\Lambda= n ! ^{(2n)^nn}}$
for the \textcolor{black}{cardinality of the} orbits of the elements of $R$, i.e., for all $u\in R$
there holds $$|\{u^k, k\in \N\}|\le\Lambda.$$
\el
\proof
As in the previous proof, we adopt~(\ref{assu}).
For $k\in\N$ let
\begin{eqnarray*}
M_k&:=& \Bigl\{\prod_{i=0}^{n-1}a_i^{\epsilon_i},\epsilon_i\in\{0,1,\ldots,k\}\Bigr\}\\
N_k&:=& \Bigl\{\sum_{\mu\in M_k}\overline{r_\mu}\,\mu,r_\mu\in\{0,1,\ldots, n ! -1\}\Bigr\}.
\end{eqnarray*}
\textcolor{black}{Observe that $|M_k|\le (k+1)^n$ and $|N_k|\le n!^{|M_k|}$.} By Proposition~\ref{fini} it follows that for $a,b\in N_k$, the sum $a+b$
also belongs to $N_k$. On the other hand, by applying~(\ref{assu}) to $u=a_j^{2}$,
$j\in\{0,1,\ldots,n-1\}$, we obtain
$$
a_j^{2n}=\sum_{i=0}^{n-1}a_i a_j^{2i},
$$
and hence, $N_{k}=N_{k-1}$ for $k\ge 2n$. 
It follows for all $u\in R$ and all $k\in\N$ that $u^k$ is of the form
$$
u^k=\sum_{i=0}^{n-1}\mu_i(k) u^j
$$
for certain coefficients $\mu_i(k)\in N_{2n-1}$ \textcolor{black}{and hence $|\{u^k, k\in\mathbb N\}| \le |N_{2n-1}|^n \le \Lambda$.}
\endproof

\bt\label{endl}
If $n:=s(R)<\infty$ then $s(R';R)\le \operatorname{lcm}(\Lambda)+\Lambda$, where \textcolor{black}{ 
$\Lambda= n ! ^{(2n)^{n}n}$}.
\et
{\bf Remarks:}

(a) Here $\operatorname{lcm}(n)$ denotes the least common multiple of the numbers in the set $\{1,2,\ldots,n\}$.

(b) Since $R'$ is contained in every subring $T$ (with $1$) of $R$,
the given bound also holds for $s(T;R)$.

{\bf Proof of Theorem~\ref{endl}}
\\
By Lemma~\ref{mul}, there exist for arbitrary $u\in R$ integers $l<k\le \color{black}{\Lambda+1}$
such that $u^k=u^l$. Thus, we have
\begin{equation*}
u^{\operatorname{lcm}(\Lambda)+\Lambda} =u^{\operatorname{lcm}(\Lambda)+\Lambda-
\frac{\operatorname{lcm}(\Lambda)}{k-l}(k-l)}=u^{\Lambda}.\tag*{$\Box$}
\end{equation*}
We conclude this section by an example of a ring $R$ which has the 
property, that $s(R)<s(R',R)$.

{\bf Example:} Let $R=\Z_2[x]/\{x^3+x^4\}$. 

The following lemma shows that for this particular ring $s(R)\le 4$.
\bl For all polynomials $P\in \Z_2[x]$ we have that
$$
x P+(1+x)P^2+P^4\equiv 0\mod (x^3+x^4).
$$
\el
\proof
We first consider the special case $P(x)=x^m$. We have to show, that
\begin{equation*}
x x^m+(1+x)x^{2m}+x^{4m}=x^{m+1}+x^{2m}+x^{2m+1}+x^{4m}\equiv  0\mod (x^3+x^4).
\end{equation*}
This is readily checked:
\begin{alignat*}{2}
m&=0:\qquad & x+1+x+1&\equiv  0\mod (x^3+x^4)\\
m&=1:  & x^2+x^2+x^3+x^4&\equiv  0\mod (x^3+x^4)\\
m&\ge 2:& x^3+x^3+x^3+x^3&\equiv  0\mod (x^3+x^4)
\end{alignat*}
Now, for arbitrary $P$, the claim follows by additivity in $\Z_2[x]$:
\begin{equation*}
x(P_1+P_2)+(1+x)(P_1+P_2)^2+(P_1+P_2)^4=
\sum_{i=1}^2 xP_i+(1+x)P_i^2+P_i^4.
\end{equation*}
\endproof
{\bf Remark:} We leave it to the reader to verify, that in fact
$s(R)=4$.

Now, we show that $s(R';R)\ge 6$.
\bl
Let $a_i\in \Z_2$ be such that  
$\sum_{i=0}^5a_ku^k=0$ in $R$ for all $u\in R$. Then $a_0=\cdots=a_5=0$.
\el
\proof
First, by choosing $u$ to be the class of $x$ in $R$ 
(which we denote by $\bar x$), we obtain
$$
a_0+a_1\bar x+a_2\bar x^2+(a_3+a_4+a_5)\bar x^3=0\qquad\text{in $R$}
$$
and hence, we conclude that $a_0=a_1=a_2=0$ and $a_3+a_4+a_5=0$.
Next, we choose $u$ to be the class of $1+x$ in $R$. Observing that
\begin{alignat*}{2}
(1+\bar x)^3&=1+\bar x+\bar x^2+\bar x^3 &&\text{ in $R$}\\
(1+\bar x)^4&=1+\bar x^4=1+\bar x^3&&\text{ in $R$}\\
(1+\bar x)^5&=1+\bar x&&\text{ in $R$}
\end{alignat*}
we have  
\begin{eqnarray*}
0&=&a_3u^3+a_4u^4+a_5u^5=\\&=&
(a_3+a_4+a_5)+(a_3+a_5)\bar x+a_3\bar x^2+(a_3+a_4)\bar x^3\text{ in $R$}
\end{eqnarray*}
which immediately implies that $a_3=a_4=a_5=0$. This completes the proof.
\endproof
Finally we prove that $s(R';R)=6$.
\bl
For all $u\in R$ it holds that $u^3+u^4+u^5+u^6=0$ in $R$.
\el
\proof
Let $u$ be the class of a polynomial $P\in\Z_2[x]$ in $R$.

First case: $P( 0)=0$. In this case, we have
\begin{alignat*}{2}
P(x)&= x Q(x)\\
P^2(x)&\equiv x^2 Q^2(x)\mod(x^3+x^4)\\
P^3(x)&\equiv x^3 Q^3(x)\equiv x^3Q( 1)\mod(x^3+x^4)\\
P^4(x)&\equiv x^4 Q^4(x)\equiv x^3Q( 1)\mod(x^3+x^4)
\end{alignat*}
and hence $P^3(x)\equiv P^4(x)$ mod $(x^3+x^4)$. This proves the claim in this case.

Second case: $P(0)=1$. In this case, we have
\begin{alignat*}{2}
P(x)&= 1+x Q( x)\\
P^2(x)&\equiv 1+ x^2 Q^2( x)\mod(x^3+x^4)\\
P^3(x)&\equiv  (1+ x Q( x))(1+ x^2 Q^2( x))\equiv \\ &\equiv
1+ xQ( x)+ x^2Q^2( x)+ x^3Q(1)\mod(x^3+x^4)\\
P^4(x)&\equiv  1+ x^4 Q^4( x)\equiv 1+ x^3Q(1)\mod(x^3+x^4)\\
P^5(x)&\equiv  (1+ xQ( x))(1+ x^3Q(1))\equiv 1+ xQ( x)\equiv P(x)\mod(x^3+x^4)
\end{alignat*}
which allows to verify the claim easily.
\endproof

\section{Two Alternative Proofs of the R\'edei-Szele Theorem}\label{alternative}
We  start with a short direct proof of Theorem~\ref{iff}. \textcolor{black}{Let $R$ be a commutative ring with unit element.} 
One implication is immediate:

Assume that $R$ is a \textcolor{black}{finite} field and $f:R\to R$. Then the \textcolor{black}{Lagrange interpolation} polynomial
$$
p(x)=\sum_{y\in R} f(y)p_y(x),
$$
where
$$
p_y(x)=\prod_{z\in R\setminus\{y\}}(x-z)\Bigl(
\prod_{z\in R\setminus\{y\}}(y-z)
\Bigr)^{-1},
$$
represents $f$. 

For the opposite implication, we assume that every function
$f:R\to R$ can be represented by a polynomial in $R[x]$.
In particular, for the function
$$
f(x):=\begin{cases} \hfill-1, &\text{ if $x=0$}\\
\hfill0,&\text{ if $x\neq 0$}\end{cases}
$$
there exists a representing polynomial 
$$\sum_{k=0}^n a_kx^k = f(x)\quad\text{ for all $x\in R$.}$$
Since $a_0=f(0)=-1$, it follows that
$$
x\underbrace{\sum_{k=1}^n a_kx^{k-1}}_{=x^{-1}} = \sum_{k=1}^n a_kx^k = 1\quad
\text{ for all $x\in R\setminus\{0\}$.}
$$
Hence, $R$ is a field. Moreover, for all $x\in R$
\begin{equation}\label{llasst}
0=xf(x)=\sum_{k=0}^n a_kx^{k+1}.
\end{equation}
The right hand side of~(\ref{llasst}) is a polynomial of degree
$n+1$ which (in the field $R$) has at most $n+1$ roots. Hence,
$|R|\le n+1$.\endproof

A second alternative proof uses the characterization of
the rings for which $s(R)=|R|$ (see Theorem~\ref{oppo}).
This condition is necessary for the property, that all
functions from $R$ to $R$ have a polynomial representative.
\textcolor{black}{In order to rule out the case $R=\Z_4$, we use the following formula from~\cite[Theorem 6, p.9]{noebi}: If $p$ is a prime number and $m\in \N$, the number of polyfunctions over $\Z_{p^m}$ is given by
$$
\Psi(p^m):=|G(\Z_{p^m})|=\exp_p\Bigl(\sum\limits_{k=1}^ms(p^k)\Bigr).
$$
Here $s$ denotes the usual number theoretic Smarandache function (see equation~\eqref{standardsmarandache}),
and $\exp_p(q):=p^q$ for better readability. 
It follows that there are $\Psi(4)=\Psi(2^2)=2^{2+4}=64$ polyfunctions over $\Z_4$, but the number of functions from $\Z_4$ to $\Z_4$ equals $4^4=256$.}
The case $R=\rho$ is ruled out by explicit verification
that $$f(x)=\begin{cases}0 & \text{  for $x\neq 0$ and}\\ 1 & \text{  for $x= 0$}\end{cases}$$ is not a
polyfunction over $\rho$\textcolor{black}{:
Since $s(\rho)=4$, it is enough to show that no polynomial $p\in\rho[x]$ of degree $\leqslant 3$ represents $f$. Suppose there is
$$
p(x) = \sum_{k=0}^3a_kx^k
$$
representing $f$. Then $p(0)=a_0=1$ and $p(a) = 1 + a_1a = 0$,
which implies that $a_1a=1$ which is impossible since $a$ does not have a multiplicative inverse.
}\endproof

\end{document}